\documentclass[12pt,a4paper]{article}  %
\usepackage{ifdraft}

\usepackage{etoolbox}
\newcommand{\bibfile}{\jobname.bib}  %
\newcommand{\universalbib}{ref.bib}
\ifdraft{\IfFileExists{\universalbib}{\renewcommand{\bibfile}{\universalbib}}{}}{}
\newcounter{cite}
\pretocmd{\cite}{\stepcounter{cite}}{}{}

\RequirePackage[mathlines]{lineno}
\ifdraft{\linenumbers}{}

\usepackage[a4paper, textwidth=16.0cm, textheight=22.0cm]{geometry}

\usepackage{amsmath,amsthm,amssymb,amsfonts}
\usepackage{mathtools}  %
\usepackage{empheq}
\usepackage{xcolor}
\usepackage[bbgreekl]{mathbbol}
\DeclareSymbolFontAlphabet{\mathbbm}{bbold}
\DeclareSymbolFontAlphabet{\mathbb}{AMSb}
\usepackage{bbm}
\usepackage{upgreek}
\usepackage{accents}
\usepackage{xspace}
\usepackage{rotating}
\usepackage{multirow,booktabs}
\usepackage[en-US]{datetime2}

\usepackage[normalem]{ulem}
\usepackage[toc,page]{appendix}

\usepackage{sectsty}
\sectionfont{\large}
\subsectionfont{\large}

\definecolor{darkblue}{rgb}{0,0.1,0.5}
\definecolor{darkgreen}{rgb}{0,0.5,0.1}
\definecolor{darkyellow}{rgb}{0.65,0.65,0.01}

\ifdraft{
    \setlength{\marginparwidth}{2.42cm}
    \usepackage[tickmarkheight=3pt,textsize=small,backgroundcolor=blue!16,linecolor=purple,bordercolor=purple]{todonotes}
}{
    \newcommand{\todo}[1]{}
    
}

\usepackage{graphicx}
\usepackage{tikz,tikzscale,pgf}
\usetikzlibrary{arrows,arrows.meta,patterns,positioning,decorations.markings,shapes}
\usepackage{pgfplots}
\usepackage{pgfplotstable}
\usepackage[justification=centering]{caption}
\usepackage{subcaption}
\usepgfplotslibrary{fillbetween}
\pgfplotsset{compat=1.11}

\ifdraft{
    \usepackage{silence}
    \WarningFilter{xcolor}{Incompatible color definition on}
    \WarningFilter{hyperref}{Draft mode on}
    \WarningFilter{refcheck}{Unused label}
    \WarningFilter{microtype}{`draft' option active}
    \WarningFilter{latex}{Writing or overwriting file} %
    \WarningFilter{latex}{Writing file} %
    \WarningFilter{latex}{Tab has} %
    \WarningFilter{latex}{Marginpar on page} %
    \WarningFilter{latex}{author given} %
}{}

\ifdraft{\usepackage{refcheck}\newcommand{\url}{\texttt}}{
    \usepackage{hyperref}
    \hypersetup{colorlinks, linkcolor=darkblue, anchorcolor=darkblue, citecolor=darkblue, urlcolor=darkblue}
    \usepackage{url}
} %
\newcommand{\email}{\texttt}

\usepackage{eqlist}
\usepackage{enumitem}
\setlist[itemize]{leftmargin=*}
\setlist[enumerate]{leftmargin=*,label=\normalfont{(\alph*)}}

\usepackage[section]{algorithm}
\usepackage{algpseudocode,algorithmicx}

\algrenewcommand\algorithmicrequire{\textbf{Input:}}
\algrenewcommand\algorithmicensure{\textbf{Output:}}
\algrenewcommand\alglinenumber[1]{\normalsize #1.}

\newtheorem{proposition}{Proposition}[section]
\newtheorem{assumption}{Assumption}[section]

\theoremstyle{definition}

\usepackage{xpatch}
\xpatchcmd{\proof}{\itshape}{\normalfont\proofnamefont}{}{}
\newcommand{\proofnamefont}{\bfseries}

\numberwithin{equation}{section}

\usepackage{microtype}
\usepackage[nobottomtitles*]{titlesec} %
\mathcode`@="8000{\catcode`\@=\active\gdef@{\mkern1mu}}

\usepackage{relsize}
\usepackage{nccmath}

\DeclareMathOperator{\bfmax}{\mathbf{max}}

\DeclareMathOperator*{\argmin}{argmin}

\DeclareMathOperator{\Span}{span}

\DeclareMathOperator{\dist}{dist}

\newcommand{\gd}{{\scriptstyle{\text{g}}}}
\newcommand{\spm}{{\scriptstyle{\text{s}}}}
\newcommand{\YY}{\mathcal{Y}}
\newcommand{\HH}{\hat{H}}
\renewcommand{\gg}{\hat{g}}
\newcommand{\ff}{\hat{f}}
\newcommand{\RR}{\mathbb{R}}

\renewcommand{\SS}{\mathcal{S}}

\newcommand{\NN}{\mathbb{N}}
\newcommand{\KK}{\mathcal{K}}

\newcommand{\sprima}{\mbox{SPRIMA}\xspace}
\newcommand{\prima}{\mbox{PRIMA}\xspace}
\newcommand{\cmt}{\mbox{CM3}\xspace}
\newcommand{\fminunc}{\mbox{fminunc}\xspace}
\newcommand{\newuoa}{\mbox{NEWUOA}\xspace}
\newcommand{\newuoas}{\mbox{NEWUOAs}\xspace}
\newcommand{\optimist}{\mbox{OptimIST}\xspace}

\newcommand{\etal}{{\mbox{et al.}}\xspace}

\xspaceaddexceptions{]\}}

\DeclareMathAlphabet{\mathsfit}{T1}{\sfdefault}{\mddefault}{\sldefault}
\SetMathAlphabet{\mathsfit}{bold}{T1}{\sfdefault}{\bfdefault}{\sldefault}

\DeclareFontFamily{U}{dutchcal}{\skewchar\font=45 }
\DeclareFontShape{U}{dutchcal}{m}{n}{<-> s*[1.0] dutchcal-r}{}
\DeclareFontShape{U}{dutchcal}{b}{n}{<-> s*[1.0] dutchcal-b}{}
\DeclareMathAlphabet{\mathlcal}{U}{dutchcal}{m}{n}
\SetMathAlphabet{\mathlcal}{bold}{U}{dutchcal}{b}{n}

\usepackage[scr=boondoxupr]{mathalfa}

\DeclareMathAlphabet{\mathpzc}{OT1}{pzc}{m}{it} %

\title{Scalable Derivative-Free Optimization Algorithms with Low-Dimensional Subspace Techniques
    \footnote{This is only a \textbf{draft containing quick notes on the main ideas} of the work.
        The final paper will be \textbf{very different}. The theory is only sketched without detailed proofs.
        The numerical results may change as the software package is still under development.}
}

\date{\today}

\author{Zaikun Zhang
    \thanks{
    Email:~\email{zaikunzhang@gmail.com}.}
}

\begin{document}

\maketitle

\begin{abstract}
    We re-introduce a derivative-free subspace optimization framework originating from Chapter 5 of the
    Ph.D.~thesis [Z. Zhang, \textit{On Derivative-Free Optimization Methods}, Ph.D.~thesis, Chinese Academy
    of Sciences, Beijing, 2012] of the author under the supervision of Ya-xiang Yuan.
    At each iteration, the framework defines a (low-dimensional) subspace based on an approximate gradient, and
    then solves a subproblem in this subspace to generate a new iterate.
    We sketch the global convergence and worst-case complexity analysis of the framework,
    elaborate on its implementation, and present some numerical results on solving problems with
    dimension as high as $10^4$ using only inaccurate function values.
\end{abstract}

\section{Introduction}
\label{sec:introduction}

Consider the unconstrained problem
\begin{equation*}
    \label{eq:opt}
    \min\{f(x) \mathrel{:} x\in \RR^n\},
\end{equation*}
where~$f\mathrel{:} \RR^n\to\RR$ is a smooth yet possibly nonconvex function.

We make the following assumption on~$f$ and will not repeat it in the sequel.
\begin{assumption}
    \label{asm:f}
    The function~$f\mathrel{:} \RR\to \RR$ is bounded from below and differentiable, and its gradient~$\nabla f$ is
    Lipschitz continuous on~$\RR^n$ with a Lipschitz constant~$L\in(0, \infty)$.
\end{assumption}

We define~$f_* =  \inf\{f(x) \mathrel{:} x\in \RR^n\}$.
For any sequence~$\{x_k\}\subset \RR^n$, we denote~$f_k = f(x_k)$ and~$g_k = \nabla f(x_k)$ for each~$k$.

We focus on derivative-free optimization (DFO) algorithms for problem~\eqref{eq:opt}.
Our target is to solve this problem with~$n$ as large as~$10^4$.
We rely on low-dimensional subspace techniques to be detailed in the sequel.

\section{A subspace framework for optimization}\label{sec:optimist}

Algorithm~\ref{alg:optimist} presents the iterated-subspace optimization framework by
Conn~\etal~\cite{Conn_Toint_Sartenaer_Gould_1996}.
Step~\ref{stp:ss} of the algorithm chooses a subspace~$\SS_k$ to explore.
When~$\dim(\SS_k)\equiv 1$, the algorithm reduces to a line search method.
Step~\ref{stp:subp} sets~$x_{k+1}$ to an approximate solution to the subspace subproblem,
which is generally much easier to tackle than the original full-space problem due to the low dimensionality.
Despite writing~``$\argmin$'', we do {not} require this subproblem to have a unique global
minimizer.%

\begin{algorithm}[htbp!]
    \caption{\label{alg:optimist} \textbf{Optim}ization with \textbf{I}terated-\textbf{S}ubspace
    \textbf{T}echnique (\textbf{\optimist})}
  {
      Input:~\mbox{$f\mathrel{:} \RR^n \to \RR$, $x_0\in\RR^n$.}
   \\For $k=0,1,2,\dots$, iterate the following.
  }
    \begin{algorithmic}[1]
        \State\label{stp:ss}Choose a subspace~$\SS_k \subset\RR^n$.%
        \State\label{stp:subp}Calculate~$x_{k+1} \approx \argmin \{f(x) \mathrel{:} x\in x_k+\SS_k\}$.
    \end{algorithmic}
\end{algorithm}

Algorithm~\ref{alg:optimist} is conceptual.
For its global convergence, we need to impose some conditions on the
subspaces~$\{\SS_k\}$ and the subproblem solutions~$\{x_{k+1}\}$. Proposition~\ref{prop:conv}
provides a necessary and sufficient condition when the objective function is convex with bounded level sets.

\begin{proposition}
    \label{prop:conv}
    Suppose that~$f$ is convex with bounded level sets.
    If Algorithm~\ref{alg:optimist} ensures~$f_{k+1} \le f_k$ for each~$k\ge 0$,
    then~$f_k \to f_*$ if and only if
\begin{equation}
    \label{eq:cond}
    \dist(g_k, \SS_k) \;\to\; 0
    \quad \text{ and } \quad
    f_{k+1} - \inf\{f(x)\mathrel{:}x\in x_k +\SS_k\} \;\to\; 0,
    \quad k \to \infty.
\end{equation}
\end{proposition}

\begin{proof}
    (a) The ``if'' part. Let~$P_k$ be the orthogonal projection onto~$\SS_k$. Then we have
    \begin{equation*}
        \begin{split}
        \frac{1}{2L}\|P_kg_k\|^2 & \;\le\; f_k - \inf\{f(x)\mathrel{:}x\in x_k +\SS_k\} \\
        & \;=\; (f_k - f_{k+1}) + \left[f_{k+1} - \inf\{f(x)\mathrel{:}x\in x_k +\SS_k\}\right]
        \;\to\; 0,
        \end{split}
    \end{equation*}
    where~$f_k-f_{k+1}\to 0$ due to the monotonicity and boundedness of~$\{f_k\}$.
    Therefore, $\|g_k\| \le \dist(g_k, \SS_k) + \|P_kg_k\| \to 0$. This implies that~$f_k\to f_*$
    since~$f$ is convex with bounded level sets.

    (b) The ``only if'' part. If~$f_k\to f_*$, then~$\|g_k\|\to 0$ by the Lipschitz continuity
    of~$\nabla f$.
    Hence~$\dist(g_k, \SS_k) \le \dist(g_k, 0) \to 0$. Meanwhile,
    \begin{equation*}
        f_{k+1} - \inf\{f(x)\mathrel{:}x\in x_k +\SS_k\} \;\le\; f_{k+1} - f_* \to 0.
    \qedhere
    \end{equation*}
\end{proof}

According to the ``if'' part of the above proof,
even without the convexity or lower-boundedness of~$f$,
condition~\eqref{eq:cond} still renders~$\|g_k\|\to 0$; if this condition holds only for~$k$ in
an infinite subset of~$\NN$, then we will have~$\liminf_k\|g_k\| = 0$ instead.

Although Proposition~\ref{prop:conv} is mainly of theoretical interest,
it suggests a general strategy to implement Algorithm~\ref{alg:optimist}:
\begin{enumerate}
    \item choose a subspace~$\SS_k$ that contains a vector~$\gg_k\approx g_k$;
    \item set~$x_{k+1}$ to a sufficiently accurate solution to~$\min \{f(x)\mathrel{:}x\in x_k +\SS_k\}$.
\end{enumerate}
The construction of~$\SS_k$ essentially demands an approximate gradient at~$x_k$. The calculation
of~$x_{k+1}$ needs more elaboration, especially if the problem is nonconvex. If we can implement
this strategy without using derivatives, then we will get a~DFO algorithm. This will be the focus of
Section~\ref{sec:dfoptimist}.

\section{A subspace framework for derivative-free optimization}
\label{sec:dfoptimist}

Algorithm~\ref{alg:dfoptimist} is a specialization of Algorithm~\ref{alg:optimist}.
As we will see, it can be implemented without using derivatives,
providing the framework for a class of subspace~DFO algorithms.%

\begin{algorithm}[htbp!]
\caption{\label{alg:dfoptimist} Derivative-free \optimist}%
  {
      Input:~\mbox{$f\mathrel{:} \RR^n \to \RR$, $x_0\in\RR^n$, $\delta_0 \in (0,\infty)$,
      $\eta \in (0,\infty)$.%
  }
   \\For $k=0,1,2,\dots$, iterate the following.
  }
    \begin{algorithmic}[1]
        \State\label{stp:g}Generate an approximate gradient~$\gg_k$ for~$f$ at~$x_k$.
        \State\label{stp:space}Choose a subspace~$\SS_k \subset\RR^n$ with $\gg_k\in \SS_k$.
        \State \label{stp:spopt}Calculate~$x_{k+1} \approx \argmin \{f(x) \mathrel{:} x\in x_k+\SS_k\}$ so
        that~$f_{k+1} \le f_k$ and%
        \begin{equation}
            \label{eq:sd}
            f_{k+1} \;\le\; {\bfmax}\left\{ f_k - \eta\delta_k^2,\;\, f(x_k - \delta_k\gg_k/\|\gg_k\|) \right\}.
        \end{equation}
        \State\label{stp:delta}If~$\|\gg_k\| \ge \eta \delta_k$ and~$f_{k+1} \le f_k -\eta \delta_k^2$,
        then~$\delta_{k+1} =  2\delta_k$\,;
        otherwise, $\delta_{k+1} = \delta_k/2$.
    \end{algorithmic}
\end{algorithm}

If~$\gg_k = 0$, then we define~$\gg_k/\|\gg_k\| = 0$ in Step~\ref{stp:spopt} of Algorithm~\ref{alg:dfoptimist}.
The adjustment of~$\delta_k$ in Step~\ref{stp:delta} resembles that of the trust region radius
in~\cite{Bandeira_Scheinberg_Vicente_2014,Gratton_Royer_Vicente_Zhang_2018}.~It can be generalized
to update~$\delta_k$ in a more refined fashion,
but this simple scheme is sufficient for our discussion.%

The \textbf{key idea} of Algorithm~\ref{alg:dfoptimist} is simple:~use the local information
to discover a subspace and then explore it.~In contrast, model-based trust-region methods
use the local information to generate a single trial point in the full space and then move on, which
may not exploit the information sufficiently.

\subsection{Global convergence and worst-case complexity}
\label{ssec:theory}

We sketch the theoretical analysis of Algorithm~\ref{alg:dfoptimist}.
To this end, we propose the following assumption on the approximate gradients~$\{\gg_k\}$.
We will elaborate on how to guarantee this assumption in Subsection~\ref{ssec:subspace}
\begin{assumption}
    \label{asm:g}
There exists a
constant~$\zeta>0$ such that
\begin{equation}
    \label{eq:g}
    \|\gg_k - g_k\|  \;\le\;  \zeta \delta_k
    \quad \text{ for each } \quad k\ge 0.
\end{equation}
\end{assumption}

Now we explain why Algorithm~\ref{alg:dfoptimist} should be globally convergent with provable
worst-case complexity bounds. We will use our proof techniques
in~\cite{Gratton_Royer_Vicente_Zhang_2018,Gratton_Vicente_Zhang_2021}.
Define
\begin{equation}
    \label{eq:KK}
    \KK \;=\;  \left\{k\in\NN \mathrel{:}
    \|\gg_k\| \ge \eta \delta_k
    \text{\; and \;}
    f_{k+1} \le f_k - \eta\delta_k^2
\right\}.
\end{equation}
Let~$\mu = 2/(L + 2\eta + 4\zeta)$.
By the triangle inequality and Taylor expansion, we can establish
\begin{equation}
    \label{eq:smalld}
    \begin{split}
        \begin{cases}
            \delta_k \le \mu\|g_k\|,& \\
        \|\gg_k - g_k \| \le \zeta\delta_k &
        \end{cases}
        \Longrightarrow \quad
        \begin{cases}
            \|\gg_k\|\ge\eta\delta_k, & \\
        f(x_k - \delta_k \gg_k/\|\gg_k\|) \le f_k - \eta \delta_k^2\;. &
        \end{cases}
    \end{split}
\end{equation}
Consequently, whenever~$\delta_k \le \mu \|g_k\|$,
we will have~$k\in\KK$ according to~\eqref{eq:sd} and~\eqref{eq:g},
and hence~$\delta_{k+1} = 2\delta_k$.
Then, as in the analysis of trust-region methods, we can obtain
\begin{equation*}
    \delta_k \;\ge\; \nu \|\tilde{g}_k\| \quad \text{with} \quad \nu = \min\{\delta_0/\|g_0\|, \; \mu/2\}
    \quad \text{and} \quad \|\tilde{g}_k\| =  \min_{0\le \ell \le k} \|g_\ell \|.
\end{equation*}
Hence we establish a lower bound for the reduction achieved by the iterations in~$\KK$, namely%
\begin{equation}
    \label{eq:red}
    f(x_{k}) - f_{k+1} \;\ge\; \eta\delta_k^2 \;\ge\; \eta\nu^2\|\tilde{g}_k\|^2
    \quad \text{for each \quad} k \in \KK.
\end{equation}
Now let~$K_\epsilon = \min \{k \in \NN \mathrel{:} \|g_k\| \le \epsilon\}$
for any given~$\epsilon > 0$. Based on~\eqref{eq:red}, we can demonstrate
that~$K_\epsilon  = \mathcal{O}(\nu^{-2}\epsilon^{-2})$ in the general nonconvex case,
$K_\epsilon = \mathcal{O}(\nu^{-2}\epsilon^{-1})$ if~$f$ is convex,
and~$K_\epsilon = \mathcal{O}(\nu^{-2}|\log\epsilon|)$ if~$f$ is strongly convex,
using the techniques in~\cite{Gratton_Royer_Vicente_Zhang_2018,Gratton_Vicente_Zhang_2021}.
These bounds also imply the global convergence of Algorithm~\ref{alg:dfoptimist}.

Note that the convergence of Algorithm~\ref{alg:dfoptimist} is essentially guaranteed by~\eqref{eq:sd}
and~\eqref{eq:g}. The subspace~$\SS_k$ does not play a direct role in the convergence analysis.
However, the choice of~$\SS_k$ is crucial for the practical performance of the algorithm.

\subsection{Defining the subspace without using derivatives}
\label{ssec:subspace}

How to generate the approximate gradient~$\gg_k$ and ensure~\eqref{eq:g}?
We can construct a model~$\ff_k$ for~$f$ in~$\RR^n$ around~$x_k$,
and then take~$\gg_k = \nabla \ff_k(x_k)$.
Condition~\eqref{eq:g} holds if~$\ff_k$ is a fully linear
model~\cite{Conn_Scheinberg_Vicente_2009b} of~$f$ in a ball
centered at~$x_k$ with a radius proportional to~$\delta_k$.
Such an~$\ff_k$ can be obtained by linear or~\mbox{(underdetermined)} quadratic Lagrange
interpolation of~$f$ on a well-poised interpolation set~$\YY_k$ in this ball~\cite{Conn_Scheinberg_Vicente_2009b}.

As an illustration, consider
\begin{equation}
    \label{eq:Y}
\YY_k = \{x_k\}\cup\{x_k+\tau \delta_k e^i\mathrel{:} i = 1,..., n\},
\end{equation}
with~$e^i$ being the~$i$th coordinate vector and~$\tau$ being a positive constant.
In this case, the linear interpolation of~$f$ on~$\YY_k$ is equivalent to a forward finite difference,
and it will provide a~$\gg_k$ satisfying~\eqref{eq:g} with~$\zeta = \tau\sqrt{n}L/2$.
Note that the function evaluations over~$\YY_k$ can be done in parallel.

We may include some previously evaluated points into~$\YY_k$ as long as they do not deteriorate
the geometry of the interpolation set. However, if Algorithm~\ref{alg:dfoptimist} converges fast enough,
it is affordable to define~$\YY_k$ like in~\eqref{eq:Y} without reusing previous points at all.

Algorithm~\ref{alg:dfoptimist} has the advantage that it converges regardless of the precise
definition of~$\SS_k$ provided that~$\gg_k \in\SS_k$.
This flexibility allows us to explore and compare different possibilities of~$\SS_k$ under
a unified framework.
In general, we will choose~$\SS_k$ according to the following principles.
\begin{enumerate}
    \item $\SS_k$ should include directions along which~$f$ is likely to decrease.
    \item $\SS_k$ should have a dimension much lower than~$n$. %
\end{enumerate}

As examples, a few possible configurations of~$\SS_k$ are listed below.
\begin{enumerate}
    \item Conjugate-gradient subspace:
        \begin{equation}
            \label{eq:yuan-stoer}
            \SS_k =\Span\{\gg_k, \, x_{k} - x_{k-1}\}.
        \end{equation}
This subspace is inspired by Yuan and Stoer's subspace perspective of conjugate-gradient
methods~\cite{Yuan_Stoer_1995}.
It is also studied by~\cite{Conn_Toint_Sartenaer_Gould_1996} and very recently
by~\cite{Zhang_Ge_Jiang_Ye_2022} in the gradient-based case.
\item Limited memory quasi-Newton
    subspace:
    \begin{equation}
        \label{eq:lmqn}
    \SS_k\;=\;\Span\{\gg_k,\, y_{k-1},\, ..., \,y_{k-m},\, s_{k-1},\, ...,\, s_{k-m}\},
    \end{equation}
where~$y_{\ell} = \gg_{\ell+1} - \gg_\ell$, $s_\ell = x_{\ell+1} - x_\ell$, and~$m\ge 1$ is an
integer much less than~$n$.
When~$\gg_k = g_k$,
this subspace is discussed by
Yuan~\cite{Yuan_2014}, and it is motivated by the fact that limited memory quasi-Newton
methods~\cite{Liu_Nocedal_1989} produce a step in this subspace.

\end{enumerate}

Furthermore, we can augment the subspaces by adding directions reflecting second-order information
of the problem if we can obtain an approximate Hessian~$\HH_k$
without using derivatives.
For instance, we can take~$\HH_k = \nabla^2\ff_k(x_k)$,~where~$\ff_k$
is the model that generates~$\gg_k$, provided that~$\ff_k$ is nonlinear.
Alternatively, we can establish~$\HH_k$ or~$\HH_k^{-1}$ by quasi-Newton formulae using the
vectors~$\{y_\ell, s_\ell\}$ mentioned earlier.
We can also explore the methods of Hessian approximation proposed by
Hare~\etal~in~\cite{Hare_JarryBolduc_Planiden_2020_grf}.
Given an approximate Hessian~$\HH_k$, we can include the following directions into~$\SS_k$.
\begin{enumerate}
    \item Approximate Newton direction, namely a direction~$d$ such that~$\HH_k d \approx -\gg_k$.
    \item Approximate negative-curvature directions.
        It is known that algorithms can benefit from exploring negative-curvature directions
        even based on inexact information~\cite{Fasano_Lucidi_2009,Hare_Royer_2022,Royer_Wright_2018}.
        If~$\HH_k$ is not positive semidefinite, we can take eigenvectors
        of~$\HH_k$ associated with negative eigenvalues as approximate negative-curvature directions.
\end{enumerate}

Recall that the conjugate-gradient method is known to be inefficient for ill-conditioned problems.
If we consider the conjugate-gradient subspace~\eqref{eq:yuan-stoer},
it is particularly important to augment the space to include the approximate Newton direction.
Otherwise, Algorithm~\ref{alg:dfoptimist} will converge slowly when the problem is ill-conditioned.

\subsection{Solving the subspace subproblem without using derivatives}
\label{ssec:subproblem}

How to calculate~$x_{k+1} \approx \argmin \{f(x) \mathrel{:} x\in x_k+\SS_k\}$ and ensure~\eqref{eq:sd}?
We propose the following strategy.
\begin{enumerate}
    \item
Invoke a DFO solver to solve approximately the low-dimensional subproblem
        \begin{equation}
            \label{eq:sub}
            \min\{f(x) \mathrel{:} x\in x_k +\SS_k\},
        \end{equation}
obtaining an approximate solution~$x_k^\spm$.
\item If~$f(x_k^\spm) \le f_k - \eta \delta_k^2$, then set~$x_{k+1} = x_k^\spm$; otherwise,
    evaluate~$f$ at~$x_{k}^\gd = x_k - \delta_k \gg_k/\|\gg_k\|$, and set~$x_{k+1}$
    to the point
    with the smallest function value
    in $\{x_k,\,x_k^\spm,\, x_k^\gd\}$.
\end{enumerate}

Note that our strategy does not impose any requirement on the quality of~$x_k^\spm$.
Indeed, the purpose of~$x_k^\spm$ is to explore the subspace~$\SS_k$, and that of~$x_k^\gd$ is to
provide a \textbf{safeguard} when the exploration fails.
It is also worth mentioning that problem~\eqref{eq:sub} is an unconstrained derivative-free optimization
problem and can be handled by, for example, \newuoa~\cite{Powell_2006}.
Take~$\SS_{k} = \Span\{\gg_k,\,x_k - x_{k-1}\}$ as an example. Problem~\eqref{eq:sub} is
equivalent to
\begin{equation*}
\min\{f(x_k + \alpha \gg_k +\beta(x_k - x_{k-1}) \mathrel{:} \alpha, \beta \in\RR\},
\end{equation*}
a low-dimensional unconstrained problem.

\subsection{Evaluation complexity and scalability}
\label{ssec:scalability}

We now examine the complexity of Algorithm~\ref{alg:dfoptimist}
in terms of function evaluations, paying particular attention to its dependence on~$n$.~Take the
nonconvex case as an example, and suppose that
the interpolation set~$\YY_k$ is $\{x_k\}\cup\{x_k+\tau \delta_k e^i\mathrel{:} i = 1,..., n\}$.~As
speculated above, the worst-case iteration complexity is~$K_\epsilon = \mathcal{O}(\nu^{-2}\epsilon^{-2})$.
If~$\tau = \mathcal{O}(n^{-p})$, then~$\zeta = \mathcal{O}(n^{\frac{1}{2}-p})$, and
hence~$\nu^{-1}= \mathcal{O}({\mu^{-1}}) =\mathcal{O}(L+2\eta+4\zeta) =\mathcal{O}(\max\{1,\, n^{\frac{1}{2}-p}\})$.
Thus $K_\epsilon = \mathcal{O}(n^{1-2p} \epsilon^{-2})$
when~$p \in [0,1/2]$.
Since~\eqref{eq:sub} is low dimensional and we have no requirement on~$x_k^{\spm}$,
we can allocate~$\mathcal{O}(1)$ function evaluations to Step~\ref{stp:spopt} for
solving~\eqref{eq:sub}. Then the complexity of function evaluations will
be~$K_\epsilon^f = \mathcal{O}(nK_\epsilon) = \mathcal{O}(n^{2-2p}\epsilon^{-2})$. In particular,
$K_\epsilon^f = \mathcal{O}(n\epsilon^{-2})$ if we take~$p = 1/2$. This linear dependence on~$n$
ensures the scalability of Algorithm~\ref{alg:dfoptimist} for the targeted problems, where~$n$
is in the order of~$10^3\sim10^4$. In contrast, the evaluation complexity of standard full-space trust-region
DFO methods is~$\mathcal{O}(n^2\epsilon^{-2})$, which is also the case for the~RSDFO method of
Cartis and Roberts~\cite{Cartis_Roberts_2022}.

\section{\sprima: A package based on \optimist and \prima}\label{sec:sprima}

\subsection{Historical remarks}

The thesis~\cite{Zhang_2012} implemented
Algorithm~\ref{alg:dfoptimist} by solving the subproblem~\eqref{eq:sub} using Powell's \newuoa~\cite{Powell_2006},
leading to the \newuoas method.
The MATLAB implementation of \newuoas was ported to Module-3 in 2016
by Dr. M. Nystr\"om (Principle Engineer at the Intel Corporation)
and made available in the open-source package~\cmt.\footnote{
\url{https://github.com/modula3/cm3/blob/master/caltech-other/newuoa/src/NewUOAs.m3}
}
It has been used by Intel in the design of chips, including its flagship product Atom P5900.
It also leads to the BBGP-sDFO method for high-dimensional analog circuit synthesis~\cite{Gu_etal_2024}.

We will extend~\newuoas to develop a new package named~\sprima, where the
subproblem~\eqref{eq:sub} is solved by the solvers in \prima~\cite{Zhang_2023}.

\subsection{Numerical experiments}

Here we present some numerical experiments based on \newuoas.\footnote{
    \url{https://github.com/newuoas/newuoas}
}

\subsubsection{Comparing \newuoas and \newuoa on moderate-dimensional problems}

\begin{figure}[htbp!]
    \centering
    \captionsetup{justification=centering}

    \begin{subfigure}[b]{0.47\textwidth}
        \includegraphics[width=\textwidth]{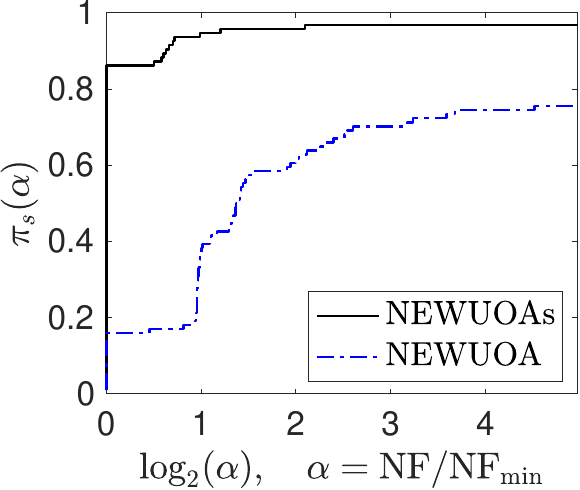}
    \vskip0.6ex
\caption{Tolerance of convergence test:~$10^{-1}$}
    \end{subfigure}
    \hfill  %
    \begin{subfigure}[b]{0.47\textwidth}
        \includegraphics[width=\textwidth]{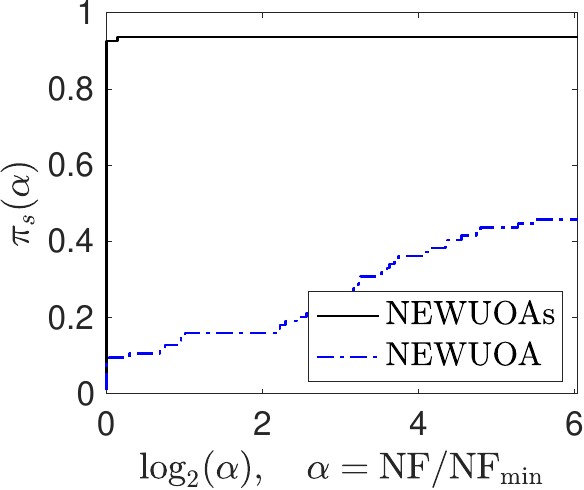}
    \vskip0.6ex
\caption{Tolerance of convergence test:~$10^{-3}$}
    \end{subfigure}
    Horizontal: relative cost in $\log_2$-scale; vertical: percentage of problems
    solved\\[0.6ex]
    \caption{
  \label{fig:newuoas}
  Performance Profiles of \newuoas and \newuoa
  \\(function values were truncated to $3$ significant digits; dimension~$n = 200$)}
\end{figure}

We tested \newuoas and \newuoa on 98 unconstrained CUTEst~\cite{Gould_Orban_Toint_2015} problems
listed in Table~\ref{tab:prob}. All the problems have changeable dimensions,
which was to~$n = 200$ in this experiment.
We provided the solvers with only the first~3 significant digits
of the function values, intending to test the robustness of the algorithms.
Figure~\ref{fig:newuoas} shows the performance profiles~\cite{Dolan_More_2002,More_Wild_2009}
generated with the tolerance of the convergence test set to~$10^{-1}$ and~$10^{-3}$.
It is evident that~\newuoas outperformed~\newuoa in this experiment.

\begin{table}[htb!]
    \captionsetup{justification=centering}
    \centering
    \begin{tabular}{llllll}
        \toprule
arglina  & arglina4  & arglinb  & arglinc  & argtrig  & arwhead \\
bdqrtic  & bdqrticp  & bdvalue  & biggsb1  & brownal  & broydn3d \\
broydn7d  & brybnd  & chainwoo  & chebquad  & chnrosnb  & chpowellb  \\
chpowells  & chrosen  & cosine  & cragglvy  & cube  & curly10  \\
curly20  & curly30  & dixmaane  & dixmaanf  & dixmaang  & dixmaanh  \\
dixmaani  & dixmaanj  & dixmaank  & dixmaanl  & dixmaanm  & dixmaann  \\
dixmaano  & dixmaanp  & dqrtic  & edensch  & eg2  & engval1  \\
errinros  & expsum  & extrosnb  & exttet  & firose  & fletcbv2  \\
fletcbv3  & fletchcr  & fminsrf2  & freuroth  & genbrown  & genhumps  \\
genrose  & indef  & integreq  & liarwhd  & lilifun3  & lilifun4  \\
morebv  & morebvl  & ncb20  & ncb20b  & noncvxu2  & noncvxun  \\
nondia  & nondquar  & penalty1  & penalty2  & penalty3  & penalty3p  \\
powellsg  & power  & rosenbrock  & sbrybnd  & sbrybndl  & schmvett  \\
scosine  & scosinel  & serose  & sinquad  & sparsine  & sparsqur  \\
sphrpts  & spmsrtls  & srosenbr  & stmod  & tointgss  & tointtrig  \\
tquartic  & tridia  & trigsabs  & trigssqs  & trirose1  & trirose2  \\
vardim  & woods  \\
\bottomrule
    \end{tabular}
    \caption{98 CUTEst problems with changeable dimensions}
    \label{tab:prob}
\end{table}

\subsubsection{Comparing \newuoas and \fminunc on large-dimensional problems}

To illustrate the scalability of \newuoas, we compared it with \fminunc on 12 problems
with~$10^4$ variables. The latter is a solver provided by MATLAB's Optimization Toolbox,
and it uses a BFGS method based on finite-difference gradients when derivatives are unavailable.
We did not test \newuoa in this case because it cannot scale to such high dimensions.
In this experiment, we still provided the solvers with only the first~3 significant digits.
The results are summarized in Table~\ref{tab:newuoas}.
Indeed, \fminunc always stopped prematurely on these problems, while~\newuoas provided reasonable
solutions.

Although we have no intention to claim that \newuoas can solve~$10^4$-dimensional problems in general,
the results in Table~\ref{tab:newuoas} demonstrate its potential for high-dimensional problems.
This illustrates the power of the subspace strategy in scaling up DFO methods.

\begin{table}[htb!]
    \captionsetup{justification=centering}
    \centering
    \caption{Comparison between \newuoas and \fminunc on $10^4$-dimensional problems\\
    \!\!\!\!\!\mbox{(function values were truncated to $3$ significant digits
    )}\\[2ex]
    }
    \label{tab:newuoas}
    \begin{tabular}{rrrrrrrrr}
        \toprule
        & && & \multicolumn{2}{c}{\newuoas} & & \multicolumn{2}{c}
        {\fminunc}\\
           \cline{5-6}  \cline{8-9}
           Problem &&  \multicolumn{1}{c}{$f(x_0)$}  &&  \multicolumn{1}{c}{$f(x_\text{fin})$}
           & \multicolumn{1}{c}{NF} && \multicolumn{1}{c}{$f(x_\text{fin})$} & \multicolumn{1}{c}{NF} \\
\hline
arwhead&&   2.99E$+$04&&  0.00E$+$00& 90331 && 2.99E+04 & 170017\\
brybnd&&    3.60E$+$05&&  4.50E$-$15& 370895 && 1.61E+04&270027\\
chrosen&&  1.99E$+$05&&  8.80E$-$14& 851736&& 1.99E+05 & 140014\\
cragglvy &&  5.49E+06 && 3.40E+03 & 110483 && 5.49E+06 & 140014 \\
dixmaane&&  7.36E$+$04&&  1.02E$+$00& 170658&& 1.23E+04 & 90018\\
engval1&&  5.89E$+$05&&  1.10E$+$04& 230880 && 5.89E+05 & 140014\\
eg2&&    8.41E$+$03&&  $-$9.99E$+$03& 110353 && $-$9.62E+03 & 160016\\
liarwhd&&  5.85E$+$06&&  7.89E$-$14& 130464 && 3.71E+05 & 240024\\
nondia&&   1.01E$+$08&&  1.97E$+$00& 90242 && 1.01E+08 & 90009\\
power&&    3.33E+11 && 1.64E$+$06 & 270951 && 3.33E+11 & 20002 \\
sparsqur&&  1.40E+07&& 1.12E$-$18& 410989 && 3.98E+06& 90009\\
woods&&  4.79E+07&& 1.97E+04& 90339 && 4.79E+07& 120012\\
\bottomrule
    \end{tabular}
    \vspace{2ex}
    \begin{center}
    \mbox{$x_0$: starting point;~\;$x_\text{fin}$: final iterate;~\;NF: number of function
    evaluations}\\
    \end{center}
\end{table}



%
%
%
\ifnum\value{cite}>0
    \small
    \bibliography{\bibfile}
    \bibliographystyle{plain}
\fi

\end{document}